\documentclass[11pt]{article}

\usepackage{latexsym,epsf}
\usepackage[dvips]{graphicx}
\usepackage{a4wide,amsthm,amssymb,amsfonts,epsfig,color,amsmath,bbm}
\usepackage{array}
\input{psfig.sty}

\theoremstyle{plain}
\newtheorem{theorem}{Theorem}
\newtheorem{lemma}[theorem]{Lemma}

\theoremstyle{remark}
\newtheorem{remark}{Remark}

\newtheorem{observation}[remark]{Observation}

\def\lcr{\operatorname{\overline{cr}}}

\long\def\comm#1{\ignorespaces}
\def\comments{\long\def\comm##1{\message{COMMENT: ##1}{\bf(( ##1 ))}}}
\comments

\ifnum\catcode`\@=\active \def\affe{@@} \else \def\affe{@} \fi

\begin{document}

\title{New lower bounds for the number of~$(\leq k)$-edges and  the rectilinear crossing number of~$K_n$
\thanks{Research on this paper was started while the first author
was a visiting professor at the Departamento de Matem\'aticas,
Universidad de Alcal\'a, Spain.}\\
} \vspace{-0.3cm}

\author{
      \normalsize Oswin Aichholzer
      \footnote{Research partially supported by the FWF
      (Austrian Fonds zur F\"orderung der Wissenschaftlichen Forschung)
      under grant S09205-N12, FSP Industrial Geometry}\\
      \small\sf Institute for Softwaretechnology\\
      \small\sf Graz University of Technology\\[-1mm]
      \small\sf Graz, Austria\\[-1mm]
      \small\sf {\tt oaich\affe ist.tugraz.at}
\and
      \normalsize Jes\'us Garc\'{\i}a
      \footnote{Research partially supported by
                grant MCYT TIC2002-01541.}\\
      \small\sf Escuela Universitaria de Inform\'atica\\
      \small\sf Universidad Polit\'ecnica de Madrid\\[-1mm]
      \small\sf Madrid, Spain\\[-1mm]
      \small\sf {\tt jglopez\affe eui.upm.es}
\and
      \normalsize David Orden
      \footnote{Research partially supported by grants MTM2005-08618-C02-02 and S-0505/DPI/0235-02.}\\
      \normalsize Pedro Ramos
      \footnote{Research partially supported by grants TIC2003-08933-C02-01 and S-0505/DPI/0235-02.}\\
      \small\sf Departamento de Matem\'aticas\\
      \small\sf Universidad de Alcal\'a \\[-1mm]
      \small\sf Alcal\'a de Henares, Spain\\[-1mm]
      \small\sf {\tt [david.orden|pedro.ramos]\affe uah.es}
}
\maketitle

\begin{abstract}
We provide a new lower bound on the number of~$(\leq k)$-edges of a
set of~$n$ points in the plane in general position. We show that
for $0 \leq k \leq \lfloor\frac{n-2}{2}\rfloor$
the number of~$(\leq k)$-edges
is at least
$$
E_k(S) \geq 3\binom{k+2}{2} + \sum_{j=\lfloor\frac{n}{3}\rfloor}^k
(3j-n+3),
$$
which, for $k\geq \lfloor\tfrac{n}{3}\rfloor$, improves the previous best
lower bound in \cite{bs-ibnks-05}.

As a main consequence, we obtain a new lower bound on the
rectilinear crossing number of the complete graph or, in other
words, on the minimum number of convex quadrilaterals determined
by~$n$ points in the plane in general position. We show that the
crossing number is at least
$$
\Bigl(\frac{41}{108}+\varepsilon \Bigr) \binom{n}{4} + O(n^3) \geq 0.379631 \binom{n}{4} + O(n^3),
$$
which improves the
previous bound of~$0.37533 \binom{n}{4} + O(n^3)$ in
\cite{bs-ibnks-05} and approaches the best known upper bound
$0.38058 \binom{n}{4}$ in \cite{AK05}.

The proof is based on a result about the structure of sets attaining
the rectilinear crossing number, for which we show that the convex
hull is always a triangle.

Further implications include improved results for small values of $n$.
We extend the range of known values for the rectilinear crossing
number, namely by $\lcr(K_{19})=1318$ and $\lcr(K_{21})=2055$.
Moreover we provide improved upper bounds on the maximum number
of halving edges a point set can have.

\end{abstract}

\medskip


\section{Introduction}
\label{sec:Introduction}

Given a graph~$G$, its \emph{crossing number} is the minimum number
of edge crossings over all possible drawings of~$G$ in the plane.
Crossing number problems have both, a long history, and several
applications to discrete geometry and computer science. We refrain
from discussing crossing number problems in their generality, but
instead refer the interested reader to the early works of
Tutte~\cite{Tutte} or Erd\H{o}s and Guy~\cite{Erdos-Guy}, the recent
survey by Pach and T\'oth~\cite{Pach-Toth}, or the extensive online
bibliography by Vrt'o~\cite{Vrto}.

In 1960 Guy~\cite{Guy} started the search for the \emph{rectilinear
  crossing number} of the complete graph,~$\lcr(K_n)$, which considers
only straight-edge drawings.  The study of $\lcr(K_n)$ is commonly
agreed to be a difficult task and has attracted a lot of interest in
recent years, see e.g.~\cite{A,AAK-02a,bs-ibnks-05,bdg-trcn-03,LVWW}.
In particular, exact values of~$\lcr(K_n)$ were only known up
to~$n=17$, see~\cite{AK05}, and also the exact asymptotic behavior is
still unknown. Several relations to other structures, like for example
\mbox{$k$-sets}, have been conjectured by Jensen~\cite{J}.
Furthermore, it has been shown by Lov\'asz et al. \cite{LVWW} that if
we denote by $E_k(S)$ the number of $(\leq k)$-edges of $S$ and
by~$\lcr(S)$ the number of crossings that appear when the complete
graph is drawn on top of $S$ (equivalently, the number of convex
quadrilaterals in $S$), then
\begin{equation} \label{eq:lvww}
\lcr(S) = \sum_{k<\frac{n-2}{2}} (n-2k-3)\,E_k(S) + O(n^3).
\end{equation}

It may be surprising that until very recently no results about the
combinatorial properties of optimal sets were known.  Motivated by
this, we start our study considering structural properties of point
sets minimizing the number of crossings, that is, attaining the
rectilinear crossing number~$\lcr(K_n)$.  Relations are obtained by
using basic techniques, like e.g. continuous motion and rotational
sweeps.
In particular, in Section~\ref{sec:Minimizing} we investigate the
changes of the order type of a point set when one of its points is
moved. We define suitable moving directions which allow us to
decrease~$\lcr(K_n)$, concluding that point configurations attaining
the rectilinear crossing number have a triangular convex hull.
Independently, and using different techniques, this result has been
extended to
pseudolinear drawings by Balogh et al.~\cite{blprs-06}.

In Section 3, and using the same technique of continuous motion, we
show that when proving a lower bound for~$(\leq k)$-edges it can be
assumed that the set has a triangular convex hull. Based on this, we
give a really simple proof of the known bound $3\binom{k+2}{2}$,
tight in the range~\mbox{$k\leq \lfloor\tfrac{n}{3}\rfloor-1$}
and, finally, we obtain a new bound for~$k\geq \lfloor\tfrac{n}{3}\rfloor$ which
improves the previous best lower bound obtained by Balogh and
Salazar~\cite{bs-ibnks-05}: We show that, for
$0 \leq k < \lfloor\frac{n-2}{2}\rfloor$, the number of~\mbox{$(\leq k)$-edges} of
a set of $n$ points in the plane in general position is at least
$$
3\binom{k+2}{2} + \sum_{j=\lfloor\frac{n}{3}\rfloor}^k
(3j-n+3).
$$
According to whether $n$ is divisible by $3$ or not, for $k\geq \lfloor\frac{n}{3}\rfloor$  this bound can be written as follows:
\begin{align*}
& 3\binom{k+2}{2} + 3\binom{k-\frac{n}{3}+2}{2} \quad \text{if $\frac{n}{3}\in\mathbb{N}$} \\[2mm]
& 3\binom{k+2}{2} + \frac{1}{3}\binom{3k-n+5}{2} \quad \text{if $\frac{n}{3}\not\in\mathbb{N}$}.
\end{align*}

If we plug our new lower bound for~$(\leq k)$-edges in Equation (\ref{eq:lvww}), we get
$$
\lcr(K_n) \geq \Bigl(\frac{41}{108}+\varepsilon \Bigr) \binom{n}{4} + O(n^3) \geq 0.379631 \binom{n}{4} + O(n^3),
$$
that improves the best previous lower bound of $0.37533 \binom{n}{4} + O(n^3)$ obtained
by Balogh and Salazar~\cite{bs-ibnks-05} and approaches the best known upper bound of
$0.38058 \binom{n}{4}$ by Aichholzer and Krasser \cite{AK05}.

For small values of $n$ the rectilinear crossing number is known for
$n \leq 17$, see~\cite{AK05} and references therein. Our results imply
that some known configurations of~\cite{A} are optimal. We thus extend
the range of known values for the rectilinear crossing number by
$\lcr(K_{19})=1318$ and $\lcr(K_{21})=2055$. Moreover our results
confirm the values for smaller $n$, especially $\lcr(K_{17})=798$,
which have been numerically obtained in~\cite{AK05}.
Finally we provide improved upper bounds on the maximum number
of halving lines
that a set of $n$ points can have.

\section{Minimizing the number of rectilinear crossings}
\label{sec:Minimizing}


Let~$S = \{p_1,...,p_n\}$ be a set of~$n$ points in the plane in
general position, that is, no three points lie on a common line. It
is well known that crossing properties of edges spanned by points
from~$S$ are exactly reflected by the order type of~$S$, introduced
by Goodman and Pollack in 1983~\cite{GP83}. The {\it order type} of
$S$ is a mapping that assigns to each ordered triple~$i,j,k$ in
$\{1,...,n\}$ the orientation (either clockwise or counterclockwise)
of the point triple~$p_i,p_j,p_k$.

Consider a point~$p_1 \in S$ and move it in the plane in a
continuous way. A change in the order type of~$S$ occurs if, and
only if, the orientation of a triple of points of~$S$ is reversed
during this process. This is the case precisely if~$p_1$ crosses the
line spanned by two other points, say~$p_2$ and~$p_3$, of~$S$.
This event has been considered previously in~\cite{AAHSW,AW} in the context of studying
the change in the number of $j$-facets under continuous motion of the points, and it is called a \emph{mutation}.

\begin{figure}[ht]
\centering
\includegraphics[scale=0.4]{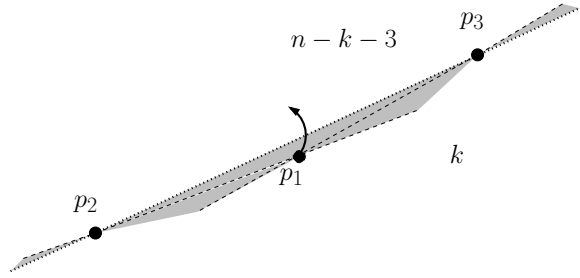}
\caption{The point~$p_1$ crosses over the segment~$p_2p_3$,
changing the orientation of the triple~$p_1,p_2,p_3$.}
\label{fig:orderflip}
\end{figure}

Assume that at time~$t_0$ the three points~$p_1,p_2,p_3$ are
collinear and that the orientation of the triple at
time~$t_0+\epsilon$ is inverse to its orientation at time
$t_0-\epsilon$ for some~$\epsilon > 0$, which can be chosen small
enough to guarantee that the orientation of the rest of triples does
not change in the interval~$[t_0-\epsilon,t_0+\epsilon]$. Let us
assume that during the mutation~$p_1$ crosses the line segment~$p_2p_3$
as indicated in Figure~\ref{fig:orderflip}; otherwise we can
interchange the role of~$p_1$ and~$p_2$ (or~$p_3$, respectively). We
say that~$p_1$ plays the \emph{center role} of the mutation.


We call the above defined mutation a \emph{$k$-mutation} if there are~$k$
points on the same side of the line through~$p_2$ and~$p_3$ as~$p_1$,
excluding~$p_1$.
Our first goal is to study how mutations affect the
number of 
crossings of~$S$, that is, the number of crossings of a
straight-line embedding of~$K_n$ on~$S$.
Note that we are considering only rectilinear
crossings.


\begin{lemma}
\label{lem:flipcross} A~$k$-mutation increases the number of crossings
of~$S$ by~$2k-n+3$.
\end{lemma}

\begin{proof}
  By definition of the~$k$-mutation, the only triple
  of points changing its orientation is~$p_1,p_2,p_3$. Thus precisely
  the~$n-3$ quadruples of points of~$S$ including this triple inverse
  their crossing properties. Observe that at time~$t_0-\epsilon$ the
  shaded region in Figure~\ref{fig:orderflip} has to be free of points
  of~$S$. Therefore, any of the~$n-k-3$ points opposite to~$p_1$ with respect to
  the line through~$p_2$ and~$p_3$ produced a crossing together with~$p_1$
  and the segment~$p_2p_3$. On the other hand, none of the~$k$ points
  on the same side as~$p_1$ does. This situation is precisely
  inverted after the flip and hence we get rid of~$n-k-3$
  crossings, but generate~$k$ new crossings.
\end{proof}


Since we know how mutations affect the number of crossings, we are now
interested in good moving directions. A point~$p \in S$ is called
\emph{extreme} if it is a vertex of the convex hull of~$S$. Two
extreme points~$p,q \in S$ are called non-consecutive if they do not
share a common edge of the convex hull of~$S$. We define a {\it
halving ray}~$\ell$ to be an oriented line passing through one
extreme point~$p \in S$, avoiding~$S \setminus \{p\}$ and splitting
$S \setminus \{p\}$ into two subsets of cardinality~$\frac{n}{2}$
and~$\frac{n-2}{2}$ for~$n$ even and~$\frac{n-1}{2}$ each for~$n$
odd, respectively. Furthermore, we orient~$\ell$ away from~$S$: For
$H$ a half plane through~$p$ containing~$S$, the `head' of~$\ell$
lies in the complement of~$H$ and the `tail' of~$\ell$ splits~$S$.

\begin{lemma}
\label{lem:rayimprove} Let~$p$ be an extreme point of~$S$ and~$\ell$
be a halving ray for~$p$. If~$p$ is moved along~$\ell$ in the given
orientation, every mutation decreases the number of crossings
of~$S$.
\end{lemma}

\begin{figure}[ht]
\centering
\includegraphics[scale=0.9]{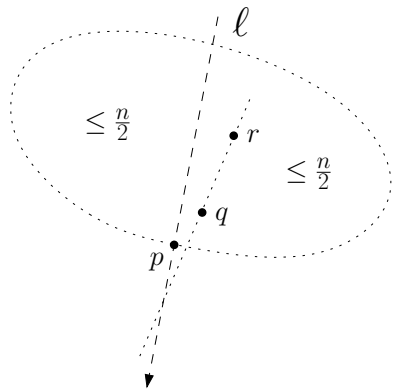}
\caption{Moving~$p$ along a halving ray~$\ell$ decreases the number
of crossings.} \label{fig:rayimprove}
\end{figure}

\begin{proof}
For the whole proof refer to Figure~\ref{fig:rayimprove}. First, we
observe that~$p$ has to be involved in any mutation and the center
role is played by another point~$q \in S$, since~$p$ is extreme. Let
$r \in S$ be the third point involved in the mutation, so that~$q$ crosses
over the segment~$pr$. As~$\ell$ is a halving ray and~$p$ an extreme
point for the~$k$-mutation which takes place when~$p$ crosses the line
defined by~$q$ and~$r$, we have that~$k \leq \frac{n}{2}-2$.
Therefore, from Lemma~\ref{lem:flipcross} it follows that the number
of crossings of~$S$ decreases.
\end{proof}


\begin{lemma}
\label{lem:halfrays} For every pair of non-consecutive extreme
points~$p$ and~$q$ of~$S$, we can choose halving rays that cross in
the interior of the convex hull of~$S$.
\end{lemma}

\begin{proof}
Let~$h$ be the line through~$p$ and~$q$. Then there is at least one
open half plane~$H$ defined by~$h$ which contains at least~$\lceil
\frac{n-2}{2} \rceil$ points of~$S$. So we can choose the two
halving rays in such a way that their tails
lie in~$H$. Now suppose that the two halving rays do not cross in
the interior of the convex hull of~$S$. Then they split~$S$ into
three regions, two outer regions and one central region. In each
outer region there are at least~$\lfloor \frac{n-1}{2} \rfloor$
points, since they are supported by halving rays. In the central
region there is at least the point which lies on the convex hull of
$S$ between~$p$ and~$q$ and not in~$H$. Finally, there are~$p$ and
$q$ themselves. All together we have~$2\cdot\lfloor \frac{n-1}{2}
\rfloor+1+2 \geq n+1$ points, a contradiction, hence the lemma
follows.
\end{proof}

\begin{observation}
Using order type preserving projective transformations it can also
be seen that a triangular convex hull can be obtained by projection
along the halving ray. This is a rather common tool when working
with order types, see e.g. \cite{Hannes}. However, we have decided
to use a self-contained, planar approach.
\end{observation}

We now have the ingredients to go for our first result, which seems
to have been a common belief (see e.g.~\cite{bdg-trcn-03}) and for
which evidence was provided by all configurations attaining
$\lcr(K_n)$ for~$n \leq 17$, \cite{A,AK05}:


\begin{theorem}
\label{thm:main2} Any set~$S$ of~$n\geq 3$ points in the plane in
general position attaining the rectilinear crossing number has
precisely~$3$ extreme points, that is, a triangular convex hull.
\end{theorem}

\begin{proof}
For the sake of a contradiction, assume that~$S$ is a set of points
attaining the rectilinear crossing number and having more than~$3$
extreme points. Let~$p$ and~$q$ be two non-consecutive extreme
points of~$S$ and let~$\ell_p$ and~$\ell_q$ be their halving rays
chosen according to Lemma~\ref{lem:halfrays}. Let~$s$ be a line
parallel to the line $h$ through~$p$ and~$q$, such that~$S$ entirely
lies on one side of~$s$ and~$\ell_p, \ell_q$ are oriented towards
$s$ (see Figure~\ref{fig:mainfig}). Furthermore, let~$s$ be placed
arbitrarily close to~$S$.

\begin{figure}[ht]
\centering
\includegraphics[scale=0.4]{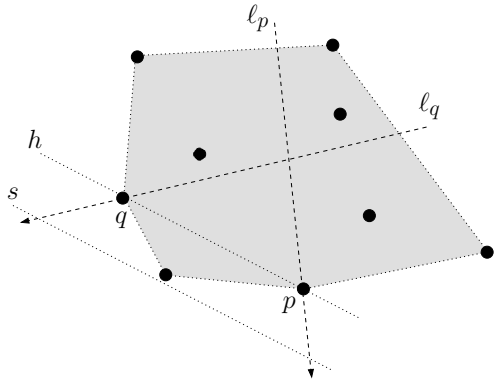}
\caption{Decreasing the number of crossings and the number of
extreme points.} \label{fig:mainfig}
\end{figure}

Now move~$p$ along~$\ell_p$
until it reaches the intersection of~$\ell_p$ and~$s$: If a mutation
occurs, we already reduce the number of crossings by
Lemma~\ref{lem:rayimprove}. Then move~$q$ along~$\ell_q$ up to the
intersection of~$\ell_q$ and~$s$. Note that, by
Lemma~\ref{lem:halfrays},~$\ell_p$ and~$\ell_q$ do not cross in
their heads; hence the movements of~$p$ along~$\ell_p$ and~$q$
along~$\ell_q$ do not interfere with each other.

After moving~$p$ and~$q$, all extreme points of~$S$ between them
changed to interior points. Hence at least one mutation happened,
since~$p$ and~$q$ are non-consecutive, and therefore the number of
crossings of~$S$ decreased, which contradicts the optimality of~$S$.
\end{proof}



\begin{observation}
If~$S$ has 3 extreme points, from the proof
of Theorem~\ref{thm:main2} it follows that we can keep moving the
three extreme points of~$S$ along their respective halving rays: If
mutations occur, this further reduces the rectilinear crossing
number. Thus, for an optimal set~$S$ the three extreme points have
to be `far away' in the following sense: For every extreme point~$p$
of~$S$, the cyclic sorted order of~$S \setminus \{p\}$ around~$p$
has to be the same as its sorted order in the direction orthogonal
to the halving ray of~$p$. (Otherwise another mutation would occur
when we keep on moving~$p$).
\end{observation}


\section{Lower bound for~$(\leq k)$-edges}
\label{s:lb}

A \emph{$j$-edge},~$0 \leq j \leq \lfloor \frac{n-2}{2} \rfloor$, is
a segment spanned by the points~$p,q \in S$ such that precisely~$j$
points of~$S$ lie in one open half space defined by the line through
$p$ and~$q$. In other words, a~$j$-edge splits~$S \setminus \{p,q\}$
into two subsets of cardinality~$j$ and~$n-2-j$, respectively.  Note
that here we consider non-oriented~$j$-edges, i.e., the edge~$pq$ equals
the edge~$qp$. We say that a $j$-edge is a \emph{halving edge} if it splits the
set as equally as possible, i.e., if $j=\tfrac{n-2}{2}$ when $n$ is even and if
$j=\tfrac{n-3}{2}$ when $n$ is odd.

A \emph{$(\leq k)$-edge} has at most~$k$ points in this half space,
that is, it is a~$j$-edge for~$0 \leq j \leq k$. We denote by
$E_k(S)$ the number of~$(\leq k)$-edges of~$S$ and omit the set when
it is clear from the context. Finally,~$(E_0,\ldots ,E_{\lfloor
\frac{n-2}{2} \rfloor})$ is the {\it~$(\leq k)$-edge vector} of~$S$.

A \emph{$k$-set} of~$S$ is a set~$S' \subset S$ of~$k$ points that
can be separated from~$S \setminus S'$ by a line (hyperplane in
general dimension). In dimension 2 there is a one-to-one relation
between the numbers of~$k$-sets and~$(k-1)$-edges, since each of
these objects can be derived from precisely two of its corresponding
counterparts. Thus, in this paper we will solely use the notion of
$j$-edges, although all the results can also be stated in terms of
$k$-sets.

In the next lemma we study how the number of~$j$-edges changes
during a mutation.

\begin{lemma}
\label{lem:flipsegment0} Let~$k\leq\tfrac{n-3}{2}$. During a
$k$-mutation, the number of~$j$-edges changes in the following way: For
$k<\tfrac{n-3}{2}$, the number of~$k$-edges decreases by one and the
number of~$(k+1)$-edges increases by one. For~$k=\tfrac{n-3}{2}$
everything remains unchanged.
\end{lemma}

\begin{proof}
  We use the same notation as for the proof of
  Lemma~\ref{lem:flipcross}. First observe that the only edges that
  change their property are the edges spanned by points~$p_1$, $p_2$, and~$p_3$. Let~$k<\tfrac{n-3}{2}$: Before the mutation, $p_1p_2$ and $p_1p_3$ are~$k$-edges,
while $p_2p_3$ is a~$(k+1)$-edge. After the mutation, the situation is reversed:
$p_1p_2$ and $p_1p_3$ are~$(k+1)$-edges while $p_2p_3$ is a~$k$-edge.
So in total we get one more~$(k+1)$-edge and one
  less~$k$-edge. For~$k=\tfrac{n-3}{2}$ the two types of
  edges considered are halving edges before and after the mutation, that is,
  the number of halving edges does not change.
\end{proof}


{From} Lemma~\ref{lem:flipcross} and
Lemma~\ref{lem:flipsegment0} we get a relation between the number
$\lcr(S)$ of rectilinear crossings of~$S$ and the number
of $j$-edges of $S$, denoted by $e_j$. An equivalent relation can
be found in~\cite{LVWW}.

\begin{lemma}
\label{lem:equation}
$$\lcr(S)+\sum_{j=0}^{\lfloor \frac{n-2}{2} \rfloor} j \cdot (n-j-2) \cdot
e_j = \frac{1}{8} \cdot (n^4-6n^3+11n^2-6n).$$
\end{lemma}

\begin{proof}
  Looking for an expression of the form $\sum_j \alpha_j e_j$ that
  cancels the variation in the number of crossings during a mutation,
  we get the relation $\alpha_{j+1}=\alpha_j+n-2j-3$. The result
  corresponds to choosing $\alpha_0=0$. The right hand side of the
  equation can be easily derived from the convex set.
\end{proof}

For the extremal case of $j$-edges, that is, halving edges, we can state a result similar to Theorem~\ref{thm:main2}:

\begin{theorem}
  For any fixed $n\geq 3$, there exist point sets with a triangular
  convex hull that maximize the number of halving edges.
\end{theorem}

\begin{proof}
As observed in the proof of Lemma~\ref{lem:rayimprove}, when an
extreme point $p$ is moved along a halving ray, only~$k$-mutations
with~$k\leq \tfrac{n}{2}-2$ can occur. Therefore, from
Lemma~\ref{lem:flipsegment0} it follows that the number of halving
lines cannot decrease. Then, we can proceed as in the proof
of Theorem~\ref{thm:main2}.
\end{proof}

One might wonder whether we can obtain a stronger result similar
to Theorem~\ref{thm:main2} stating that any point set maximizing
the number of halving edges has to have a triangular convex hull.
But there exist sets of $8$ points with $4$ extreme points bearing
the maximum of $9$ halving edges, see~\cite{AAHSW}, and similar
examples exist for larger $n$. Hence, the stated relation is tight
in this sense. We leave as an open problem the existence of a
constant $h$ such that any point set maximizing the number of
halving edges has at most $h$ extreme points. We conjecture that
such a constant exists, and the results for $n \leq 11$ suggest
that $h=4$ could be the tight bound.

Similar arguments as above can be used to prove the next result, which
is our starting point for the lower bound of~$(\leq k)$-edges:

\begin{lemma}
\label{lem:facettriangular} Let~$S$ be a set of~$n$ points with
$h>3$ extreme points and~$(\leq k)$-edge vector
$(E_0,\ldots,E_{\lfloor \frac{n-2}{2} \rfloor})$. Then there exists
a set~$S'$ of~$n$ points with triangular convex hull and
$(\leq k)$-edge vector~$(E'_0,\ldots,E'_{\lfloor \frac{n-2}{2} \rfloor})$
with~$E'_i \leq E_i$ for all~$i=0,\ldots,\lfloor \frac{n-2}{2}
\rfloor$ (where at least one inequality is strict).
\end{lemma}

\begin{proof}
  The proof follows the lines of the proof for Theorem~\ref{thm:main2}
  to obtain a set with only~$3$ extreme vertices. Observe that for
  all~$k$-mutations which occur during this process it holds~$k\leq \tfrac{n}{2}-2$
  because we are moving along halving rays. Thus by Lemma~\ref{lem:flipsegment0}
  every mutation decreases the number of
 ~$k$-edges by one and increases the number of~$(k+1)$-edges by one. For the
 ~$(\leq k)$-edge vector this means that~$E_k$ is decreased by
  one and the rest of the vector remains unchanged. The statement
  follows.
\end{proof}

As a warm-up, we start with a really simple and geometric proof of
the following bound, which has been independently shown
in~\cite{af-lbrcn-05,LVWW} using circular sequences:

\begin{theorem}
\label{thm:facets1}
  Let~$S$ be a set of~$n$ points in the plane. The number of
 ~$(\leq k)$-edges of~$S$ is at least~$3 \binom{k+2}{2}$ for~$0 \leq k < \frac{n-2}{2}$.
  This bound is tight for~$k \leq \lfloor \frac{n}{3} \rfloor -1$.
\end{theorem}

\begin{proof}
  By Lemma~\ref{lem:facettriangular} we can assume that~$S$ has a
  triangular convex hull, as otherwise we can find a point set with a
  strictly smaller~$(\leq k)$-edge vector for which the theorem still
  has to hold. Let~$p,q,r$ be the three extreme points of~$S$.

  By rotating a ray around each extreme point of~$S$, we get
  exactly three~$0$-edges and six~$j$-edges for every~$1 \leq j <
  (n-2)/2$, all of them incident to~$p$,~$q$ or~$r$. This gives a total
  of~$3+6k$~$(\leq k)$-edges, which already proves the lower bound for~$k=1$.
  For~$2 \leq k < \frac{n-2}{2}$ we will prove the lower bound by
  induction on~$n$.  The cases~$n \leq 3$ are obvious and serve as an
  induction base. So for~$n \geq 4$ consider
 ~$S_1=S\smallsetminus\{p,q,r\}$ where~$n_1=n-3 \geq 1$ denotes the
  cardinality of~$S_1$.

  Observe that, since the convex hull of~$S$ is a triangle,
  a~$j$-edge of~$S_1$ is either a~$(j+1)$-edge or a~$(j+2)$-edge of
 ~$S$. Therefore, if~$2 \leq k<\tfrac{n-2}{2}$ we get

 ~$$
  E_k(S) \geq E_{k-2}(S_1) + 3+6k \geq 3\binom{k}{2} + 3+6k = 3\binom{k+2}{2}.
 ~$$

  Finally, the example in \cite{ehss-89}
   shows that the bound~$3 \binom{j+2}{2}$ is tight
  for~$j \leq \lfloor \frac{n}{3} \rfloor -1$.
\end{proof}

In view of the preceding proof, it is clear that in order to improve
the bound for~$k\geq \lfloor \frac{n}{3} \rfloor$ we need to show
that a number of~$j$-edges of~$S_1$ are~$(j+1)$-edges of~$S$. This
is going to be our next result, but first we need some preparation.

It is more convenient now to consider {\em oriented} $j$-edges: an
oriented segment $pq$ is a $j$-edge of $S$ if there are exactly $j$
points of $S$ in the open half plane to the right of $pq$.
Following~\cite{elss}, we rotate a directed line~$\ell$ around
points of the set~$S$, counterclockwise, and in such a way that,
when~$\ell$ contains only one point, it has exactly~$k$ points
of~$S$ on its right. We refer to this movement as a
{\em~$k$-rotation}. If the line rotates around a point~$p$, the
half-lines into which~$p$ divides~$\ell$ are the {\em head} and the
{\em tail} of the line. Observe that, if during a~$k$-rotation the
line~$\ell$ reaches a new point~$q$ on its tail, then~$qp$ is
a~$(k-1)$-edge and the~$k$-rotation
continues around~$q$. On the other hand, if the new point~$q$
appears on the head of the ray, then~$pq$ is a~$k$-edge and
the~$k$-rotation continues also around~$q$. We recall that when
a~$k$-rotation of~$2 \pi$ is completed, all~$k$-edges of~$S$ have
been found.

We denote by~$\ell^+$ and~$\overline{\ell^+}$, respectively, the
open and closed half-planes to the right of~$\ell$ and, similarly,
$\ell^-$ and~$\overline{\ell^-}$ will be the half-planes to the left
of~$\ell$.

\begin{theorem} \label{pro:lb}
  Let~$S$ be a set of~$n$ points in the plane in general position and
  let~$T$ be a triangle containing~$S$. If~$\lfloor\tfrac{n}{3}\rfloor
  \leq k \leq \tfrac{n}{2}-1$, then there exist at least \hskip0.1em
 ~$3k-n+3$ \hskip0.1em~$k$-edges of~$S$ having to the right only one
  vertex of~$T$.
\end{theorem}

\begin{proof}
Let~$p$,~$q$ and~$r$ be the vertices of~$T$ in counterclockwise
order. Throughout this proof, we will refer to a~$k$-edge and its
supporting line synonymously. Moreover, edges having one vertex
of~$T$ on its right will be called {\em good} edges, and the rest
will be said to be {\em bad}.

We start with the case of halving lines for~$n$ even, which is
straightforward. There are at least~$n$ halving lines and exactly
half of them are good (because each edge is a halving line in both
orientations). Therefore, because~$k=\tfrac{n}{2}-1$, we have
that~$\tfrac{n}{2} = 3k-n+3$. In the following,~$k<\tfrac{n}{2}-1$.

Since the number of~$k$-edges is always at least~$2k+3$ (see
\cite{LVWW}), if all~$k$-edges are good the result is true.
Therefore, without loss of generality, we can assume that there are
bad edges having~$q$ and~$r$ on its right. Among them, let~$\ell_1$
be the bad~$k$-edge which intersects~$pq$ closest to~$q$.
Now, we make a~$k$-rotation of~$\ell_1$ and distinguish two cases.
For the whole proof refer to Figure~\ref{fig:triangle}.

\begin{figure}
\centering
\includegraphics{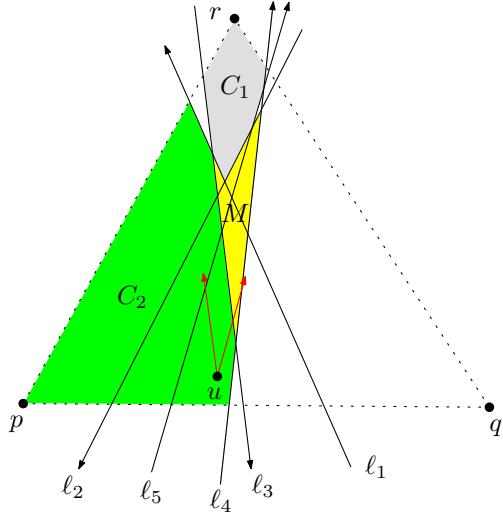}
\caption{Proving Theorem~\ref{pro:lb}:~$k$-edges and their
relation to the triangle~$T$.} \label{fig:triangle}
\end{figure}

{\it Case 1.} If a~$k$-edge having~$p$ and~$r$ on its right is not
found, then a good~$k$-edge can be found for each of the~$k$ points
to the right of~$\ell_1$: if~$a\in \ell_1^+$, consider a directed
line through~$a$ and parallel to~$\ell_1$ and rotate it around~$a$.
Before rotating 180 degrees, a~$k$-edge is found and it has to be
good because there is no~$k$-edge having~$p$ and~$r$ on its right.
Therefore, because~$k\leq\tfrac{n-3}{2}$, it holds that~$k\geq
3k-n+3$ and the result follows.

{\it Case 2.} Let~$\ell_2$ be the first~$k$-edge we obtain from the
rotation for which~$p$ and~$r$ lie on its right. Let~$H=S\cap
\ell_2^+\cap \ell_1^-$ and denote by~$h$ the cardinality of~$H$.
Observe that all~$k$-edges between~$\ell_1$ and~$\ell_2$ we get
during the rotation are good edges. Since points in~$H$ are
necessarily encountered in the head of the ray during the
$k$-rotation at least once, there is a good~$k$-edge incident to
each of them. Consider~$C_1=S\cap \overline{\ell_1^+}\cap
\overline{\ell_2^+}$ and denote by~$c_1$ its cardinality. If $c_1=0$, then
$h=k$ and the result follows as in Case 1. Let
$\ell_3$ be a~$k$-edge tangent to~$C_1$ at only one point and
leaving~$C_1$ on its left. Observe that~$\ell_3$ could have~$r$ to
its right and that it can be found by rotating a tangent to~$C_1$
counterclockwise and, if a~$k$-edge~$uv$ defined by two points of
$C_1$ is found, proceeding with the rotation with~$v$ as new center
(maybe repeatedly). Finally,
let~$C_2=S\cap\ell_1^-\cap \overline{\ell_3^+}$
and~$\ell_4$ be the common tangent to~$C_1$ and
$C_2$ leaving both sets to its left.

Now we want to bound the number of points in~$\ell_4^+$. To this
end, observe that~$k-h \leq c_1 \leq k-h+2$, depending on the number
of points defining~$\ell_2$ that belong to~$C_1$. Therefore, if we
denote by~$m$ the number of points in~$M=(S\cap\ell_4^-)\smallsetminus(C_1\cup C_2)$,
then~$|S\cap \overline{\ell_4^-}| \geq k+m+c_1+1$
(start from~$\ell_3$: it has~$k$ points to the right. In addition we
have~$m+c_1$ points plus a point not in~$C_1$ defining~$\ell_3$).
Therefore,~$|S\cap \ell_4^+|\leq n-k-m-c_1-1$. Again we have to
distinguish two cases:

{\it Case 2a.}~$c_1 \geq k-h+1$. In this case,~$|S\cap \ell_4^+|
\leq n-2k+h-m-2$. Therefore, if~$|S\cap \ell_4^+|>k$, then~$h >
3k-n+2+m$, implying~$h \geq 3k-n+3$ and thus the~$h$ good~$k$-edges
incident to points in~$H$ (see above) are sufficient to guarantee
the result.

On the other hand, if there are at most~$k$ points to the right of
$\ell_4$, we can rotate~$\ell_4$ around~$C_1$, clockwise, and find a
$k$-edge~$\ell_5$ (which could be bad and which could also coincide
with~$\ell_4$). Observe that when rotating from~$\ell_4$ to~$\ell_5$
only points from~$M$ and~$C_2$ can be passed by the line and that
one point spanning~$\ell_5$ belongs to either~$M$ or~$C_2$. Thus
$|C_2\cap \overline{\ell_5^+}| \geq k-|S\cap \ell_4^+|-m+1 \geq
3k-n+3-h$.

We finally claim that for each point in~$C_2\cap
\overline{\ell_5^+}$ we can find a good~$k$-edge, which together
with the~$h$ good edges from~$H$ settles Case 2a. To prove the claim
let~$u\in C_2\cap \overline{\ell_5^+}$ and consider the half-cone
with apex at~$u$, edges parallel to~$\ell_3$ and~$\ell_5$ and
containing~$C_1$. Because~$\ell_5$ and~$\ell_3$ are~$k$-edges, we
can guarantee that if we rotate a line around~$u$ we find at least
two~$k$-edges,~$uv$ and~$wu$, and that at least one of them is good:
Start from the line parallel to~$\ell_5$ which has less than~$k$
points to its right and more than~$k$ points to its left.  Rotating
counterclockwise around~$u$ until the line is parallel to~$\ell_3$
reverses the situation. Thus, during the rotation we first get an
edge~$uv$ with~$k$ points to its right and then an edge~$uw$
with~$k$ points to its left. If~$r$ is to the left of~$uv$ then~$uv$
is the good~$k$-edge for~$v$. Otherwise~$r$ has to be to the right
of~$uw$ and thus~$wu$ is the good~$k$-edge for~$v$. Finally observe
that all good~$k$-edges associated to points in~$C_2\cap \overline{\ell_5^+}$ 
in this last step are different from the good~$k$-edges incident to points in~$H$ 
which were found in the first part of the $k$-rotation because the former ones have $r$ 
to its left while the later ones have $r$ to its right. 

{\it Case 2b.}~$c_1 = k-h$. In this case, the arguments of Case 2a
give~$3k-n+2$ good~$k$-edges and one more is needed. Observe that, in
this case, the points defining~$\ell_2$ are to the left of~$\ell_1$.
Therefore, the point~$t$ in the tail of the~$k$-edge defining~$\ell_2$
has a good~$k$-edge incident to it: Points to the left of~$\ell_1$ and
defining~$k$-edges are found during the $k$-rotation, and the first
time they are found, they have to define a good~$k$-edge (recall
that~$\ell_2$ was the first bad edge). As $t$ does not belong
to~$H$, the good~$k$-edge incident to~$t$ (having $r$ to its right) 
was not counted previously.
\end{proof}

\begin{theorem} \label{th:lb}
  Let~$S$ be a set of~$n$ points in the plane in general position and
  let~$E_k(S)$ be the number of~$(\leq k)$-edges in~$S$.
  For $0 \leq k < \lfloor\frac{n-2}{2}\rfloor$ we have
$$
E_k(S) \geq 3\binom{k+2}{2} + \sum_{j=\lfloor\frac{n}{3}\rfloor}^k
(3j-n+3).
$$

\end{theorem}

\begin{proof}
The proof goes by induction on~$n$. Observe that
Lemma~\ref{lem:facettriangular} guarantees that it is sufficient to
prove the result for sets with triangular convex hull. Let~$p,q,r$
be the vertices of  the convex hull of~$S$ and
let~$S_1=S\smallsetminus \{p,q,r\}$.
For $k\leq \lfloor\tfrac{n}{3}\rfloor -1$ the result is already given
by Theorem~\ref{thm:facets1}.  If~$k\geq \lfloor\tfrac{n}{3}\rfloor
+1$ then
$$
E_{k-2}(S_1)\geq 3\binom{k}{2}
+ \sum_{j=\lfloor\frac{n-3}{3}\rfloor}^{k-2} (3j-(n-3)+3)=
3\binom{k}{2} + \sum_{i=\lfloor\frac{n}{3}\rfloor}^{k-1} (3i-n+3).
$$
Furthermore, as in the proof of Theorem~\ref{thm:facets1} we know
that there are exactly~$3+6k$~$(\leq k)$-edges of~$S$ adjacent
to~$p$,~$q$ and~$r$, so using Theorem~\ref{pro:lb} we conclude that
$$
E_k(S) \geq E_{k-2}(S_1) + 3 + 6k + 3(k-1)-(n-3)+3 \geq
3\binom{k+2}{2} + \sum_{j=\lfloor\frac{n}{3}\rfloor}^k (3j-n+3).
$$
For~$k = \lfloor\tfrac{n}{3}\rfloor$,~$E_{k-2}(S_1)\geq
3\binom{k}{2}$ and then
$$
E_k(S) \geq 3\binom{k}{2} + 3 + 6k + 3(k-1)-(n-3)+ 3 =
3\binom{k+2}{2} +3\Bigl\lfloor\frac{n}{3}\Bigr\rfloor-n+3.
$$
\end{proof}

As a main consequence of Theorem~\ref{th:lb}, we can obtain a new lower bound
for the rectilinear crossing number of the complete graph:

\begin{theorem}
For each positive integer $n$,
$$
\lcr(K_n)\geq\Bigl(\frac{41}{108}+\varepsilon \Bigr) \binom{n}{4} + O(n^3) \geq 0.379631 \binom{n}{4} + O(n^3).
$$
\end{theorem}

\begin{proof}
As shown in \cite{LVWW}, the number of~$(\leq k)$-edges and the crossing
number of $K_n$ are strongly related. More precisely, if we denote by~$\lcr(S)$ the number of
crossings when the complete graph is drawn with set of vertices~$S$, then
\begin{equation*}
\lcr(S) = \sum_{k<\frac{n-2}{2}} (n-2k-3)\,E_k(S) + O(n^3). \tag{1}
\end{equation*}
Writing
$3\binom{k+2}{2}+\sum_{j=\lfloor\frac{n}{3}\rfloor}^k (3j-n+3) =\hat{E}_k$,
we get
\begin{equation*}
\lcr(K_n) \geq \sum_{k<\frac{n-2}{2}} (n-2k-3)\, \hat{E}_k=
\frac{41}{108} \binom{n}{4} + O(n^3) .
\end{equation*}

Now, we can slightly improve the lower bound by exploiting a bound for
$(\leq k)$-edges which is better than $\hat{E}_k$ when $k$ is close to $\tfrac{n}{2}$: In \cite{LVWW}  it is shown that
$$
E_k(S) \geq \binom{n}{2}-n\sqrt{n^2-2n-4k(k+1)} = F_k.
$$
A straightforward computation shows that, for $n$ large enough, $F_k \geq \hat{E}_k$ if
$k\geq 0.4981 n$.

\vskip1mm

Applying again Equation~(\ref{eq:lvww}) we get
$$
\lcr(K_n)\geq \sum_{k<\frac{n-2}{2}} (n-2k-3)\, \hat{E}_k +
\sum_{k= 0.4981 n}^{\frac{n-2}{2}} (n-2k-3)\, (F_k- \hat{E}_k) +O(n^3)
= \Bigl(\frac{41}{108}+\varepsilon \Bigr) \binom{n}{4} + O(n^3).
$$
In order to give an estimation for $\varepsilon$, let $t_0=0.4981$ and observe that
\begin{equation*}
\begin{split}
\sum_{k= t_0 n}^{\frac{n-2}{2}} (n-2k-3)\, (F_k- \hat{E}_k) & =
n^3 \sum_{k= t_0 n}^{\frac{n-2}{2}} \Bigl( 1-2\,\frac{k}{n} \Bigr)  \left(\frac{1}{3}+ \frac{k}{n} -3\Bigl(\frac{k}{n}\Bigr)^2  -\sqrt {1-4\,\Bigl(\frac{k}{n}\Bigr)^2}  \right) \\
& = n^4 \, \int_{t_0}^{1/2} (1-2\,t)\,( \frac{1}{3}+t-3t^2-\sqrt{1-4\,{t}^{2}})\,dt +O(n^3).
\end{split}
\end{equation*}
Therefore,
$$
\varepsilon = 24\, \int_{t_0}^{1/2} (1-2\,t)\,\Bigl( \frac{1}{3}+t-3t^2-\sqrt{1-4\,{t}^{2}}\Bigr)\,dt \simeq
1.4 \cdot 10^{-6}.
$$

\end{proof}

\begin{observation}
  Using Theorem~\ref{th:lb} and Lemma~\ref{lem:equation} it can be
  shown that the configurations of~$19$ and $21$ points in \cite{A}
  are optimal for the number of crossings: their $(\leq k)$-edge
  vectors are, respectively, $(3,9,18,30,45,63,86,115,171)$ and
  $(3,9,18,30,45,63,84,111,144,210)$. Because they match the bound in
  Theorem~\ref{th:lb} for $k < \tfrac{n-3}{2}$, we have that
  $\lcr(K_{19})=1318$ and $\lcr(K_{21})=2055$.
\end{observation}

\begin{observation}
  Let us recall that, for $n$ odd, $j$-edges with $j=\tfrac{n-3}{2}$
  are halving edges. Let $h_n=\max_{|S|=n}
  e_{\lfloor\frac{n-2}{2}\rfloor}$ be the maximum number of halving
  lines that a set of $n$ points can have. In Table~1 we present a
  summary of the values of $h_n$ for $13\leq n \leq 27$: the value
  $h_{14}=22$ and the upper bound for $h_{16}$ were shown in
  \cite{br-02}, while the lower bound for $h_{16}$ appeared in~\cite{e-92}.
  The rest of the lower bounds come from the examples
  in~\cite{A} and the upper bounds can be derived applying
  Theorem~\ref{th:lb} with $k=\lfloor\tfrac{n-2}{2}\rfloor - 1$,
  namely $h_n \leq {n \choose 2} - E_{\lfloor\tfrac{n-2}{2}\rfloor - 1}$.
\end{observation}

\begin{table}\label{tab2}
\setlength{\extrarowheight}{1mm}
\begin{center}
\begin{tabular}{m{4mm}|m{4mm}|m{4mm}|m{4mm}|m{4mm}|m{4mm}|m{4mm}|m{4mm}|m{4mm}|
m{4mm}|m{4mm}|m{4mm}|m{4mm}|m{4mm}|m{4mm}|m{4mm}}
$n$ & 13 & 14 & 15 & 16 & 17 & 18 & 19 & 20 & 21 & 22 & 23 & 24 & 25 & 26 & 27 \\ \hline
$h_n$ & 31 & 22 & 39 & 27 \hskip3mm 28 & 47 & 33 \hskip3mm 36 & 56 & 38 \hskip3mm 43 & 66 & 44 \hskip3mm 51 & 75 \hskip3mm 76 & 51 \hskip3mm 60 & 85 \hskip3mm 87 & 57 \hskip3mm 69 & 96 \hskip3mm 99 \\
\end{tabular}
\end{center}
\caption{Values and bounds of $h_n$ for $13\leq n \leq 27$.}
\end{table}

\section{Concluding Remarks}

In this paper we have presented a new lower bound for the number of
$(\leq k)$-edges of a set of $n$ points in the plane in general
position. As a corollary of this, a new lower bound for the
rectilinear crossing number of $K_n$ is obtained. The basis of the
technique is a property about the structure of sets minimizing the
number of $(\leq k)$-edges or the rectilinear crossing number: such
sets have always a triangular convex hull.

There are still a host of open problems and conjectures about these
and related questions, among which we emphasize the following:

\begin{itemize}
\item Prove that the new lower bound is optimal for some range of $k$.
  Based on computational experiments, we conjecture that the bound is
  optimal for $k \leq \lfloor\tfrac{5n}{12}\rfloor-1$.
\item Prove that all sets maximizing the number of halving lines have
  a convex hull with at most $h$ vertices. We conjecture that $h=4$ is
  sufficient.
\item Prove that sets minimizing the crossing number maximize the
  number of halving lines.
\end{itemize}

\end{document}